\documentclass[a4paper,aip,tightenlines,preprint,onecolumn]{revtex4-1}
\usepackage[utf8]{inputenc}
\usepackage{amsmath}
\usepackage{amssymb}
\usepackage{amsthm}
\usepackage{mathtools} 
\usepackage{enumerate}
\usepackage[dvipsnames]{xcolor}
\usepackage{caption}
\usepackage[left=2.5cm,right=2.5cm,top=2cm,bottom=2cm,includeheadfoot]{geometry}

\newtheorem*{theorem*}{Theorem}

\PassOptionsToPackage{hyphens}{url} 
\usepackage[breaklinks=true]{hyperref} 
\usepackage{natbib}
\usepackage{cleveref}

\begin{document}
\title{On a combinatorial problem in the Secret Santa ritual}

\author{Alexander Steinicke}
\affiliation{Montanuniversität Leoben, Austria}

\author{Markus Penz}
\affiliation{Max Planck Institute for the Structure and Dynamics of Matter, Hamburg, Germany}
\affiliation{Basic Research Community for Physics, Leipzig, Germany}

\author{Bine Penz}
\affiliation{SalTo Vocale Chor, Salzburg, Austria}

\begin{abstract}
\begin{center}
\includegraphics[width=.82\linewidth]{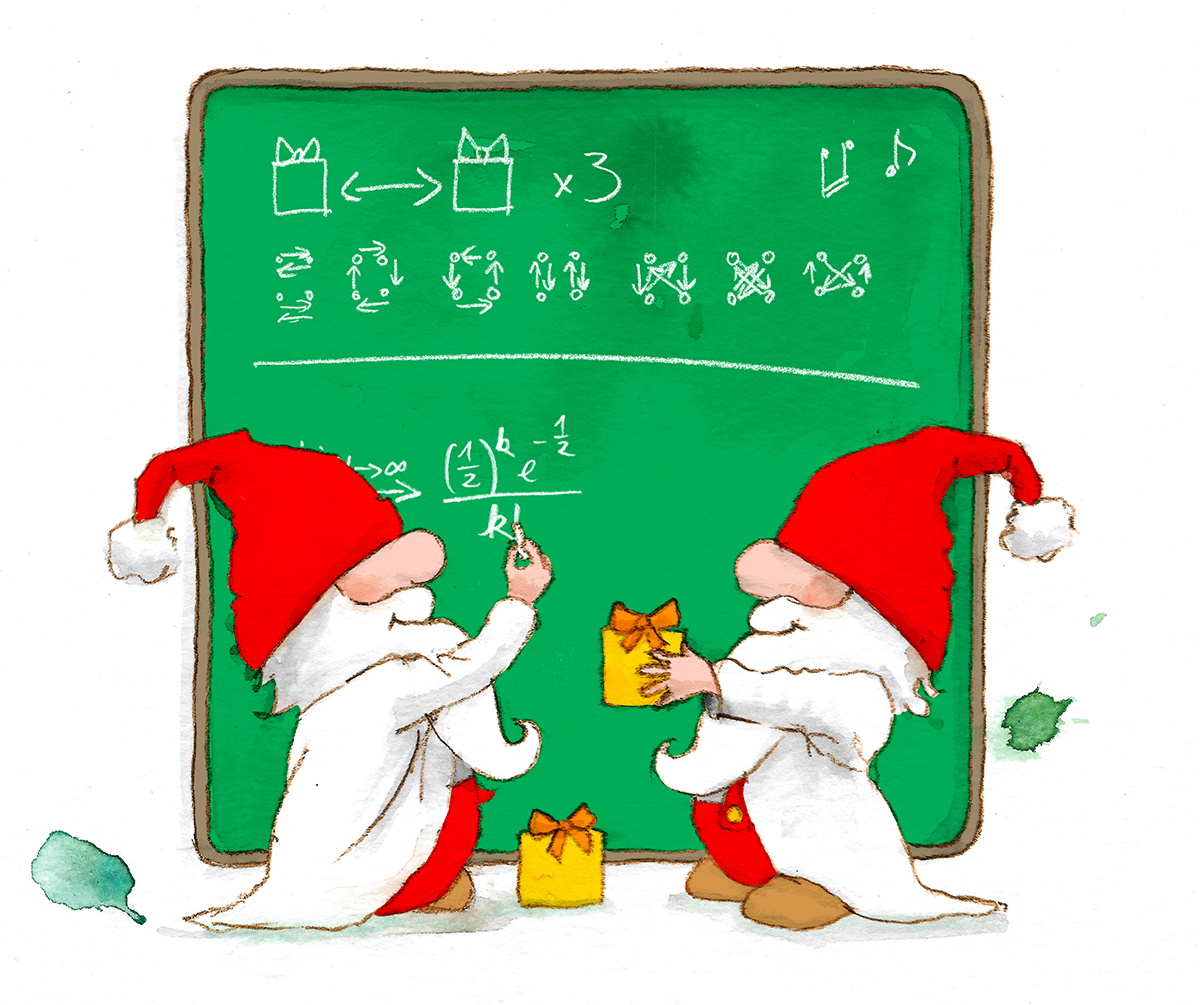}
\end{center}
\par

\noindent\emph{English.} The Secret Santa ritual, where in a group of people every member presents a gift to a randomly assigned partner, poses a combinatorial problem when considering the probabilities involved in the formation of pairs, where two persons exchange gifts mutually. We give different possible derivations for such probabilities by counting fixed-point-free permutations with certain numbers of 2-cycles.\\

\noindent\emph{German.} Das vorweihnachtliche Wichtel-Spiel, bei dem in einer Gruppe von Personen jedes Mitglied einem zufällig zugeordneten, anderen Mitglied ein Geschenk überreicht, wirft ein kombinatorisches Problem auf, wenn die Wahrscheinlichkeiten für das Auftreten von Paaren, wo zwei Personen sich gegenseitig beschenken, berück\-sichtigt werden. Wir geben verschiedene Ableitungen für solche Wahrscheinlichkeiten an, indem fixpunktfreie Permutationen mit einer bestimmten Anzahl von 2-Zyklen gezählt werden.\\

\noindent This work is dedicated to the public domain. \texttt{(CC0)}
\end{abstract}

\maketitle

\section{Introduction}

This christmas eve one of the authors put forward a combinatorial problem that originates from a gift-exchange ritual that in this case was practiced in her choir \emph{SalTo Vocale} consisting of $N=20$ persons. The ritual is commonly known as ``Secret Santa'' \cite{wiki-secret-santa}, or ``Wichteln'' in German. The version we treat here is played as follows: In a group of $N$ people everyone is randomly assigned a partner to whom at a special event a small gift is presented. Now there is of course the possibility that the person assigned to someone is just the one who has picked the other in return. Such a constellation will be referred to as a \emph{pair}. The question raised by members of the choir was then: How large is the probability that at least three pairs are formed, a situation that actually had occurred that year. In order to give a full answer to the question, we will translate it into a problem of permutation groups, then find the total number of possible configurations in a Secret Santa ritual of $N \geq 2$ persons, and finally count the number of configurations with a given number of pairs. We will see that the correct answer is given by the number of fixed-point-free permutations $F_N$ of $N$ elements and the related number of fixed-point-free permutations with a certain number of 2-cycles.

\section{Fixed-point-free permutations}\label{sec:counting-fixed-point-free}

Let the participating persons be represented by the numbers $\{1, 2, \ldots, N\}$m then the modeling of admissible configuration in a Secret Santa ritual is given by the fixed-point-free permutations (derangements), since nobody can pick themselves. For $N=2$ there is obviously only one possible configuration consisting of one pair, $1 \mapsto 2, 2 \mapsto 1$, or $(21)$ in canonical cycle notation with the largest element always in front. For $N=3$ the possible configurations are $(312)$ and $(321)$ and pairs are impossible. Now since pairs are equivalent to 2-cycles the transformation of the initial question into group theoretical/combinatorial terms is already achieved.

In counting the number of fixed-point-free permutations $F_N$ within the set of all permutations of $N$ elements $S_N$, we first define the set of permutations that include a fixed-point mapping $(j)$, $S_{N,j} \subseteq S_N$. An intersection of $k$ such sets has $k$ given fixed-points, so it will still contain $(N-k)!$ elements, and one has $N \choose k$ possibilities of choosing the fixed-points. By the inclusion-exclusion principle one thus has \cite{wiki-derangement}
\begin{equation}
|S_{N,1} \cup \cdots \cup S_{N,N}| = \sum_{k=1}^N (-1)^{k+1} {N \choose k} (N-k)! = \sum_{k=1}^N (-1)^{k+1} \frac{N!}{k!}.
\end{equation}
Now clearly the number of fixed-point-free permutations is the total number of permutations minus the above result, so
\begin{equation}\label{FN-num}
|F_N| = N! - \sum_{k=1}^N (-1)^{k+1} \frac{N!}{k!} = \sum_{k=0}^N (-1)^k \frac{N!}{k!}.
\end{equation}
Extending the sum to infinity gives a direct relation to Euler's $e$,
\begin{equation}\label{eq:FN-Euler}
|F_N| = \sum_{k=0}^\infty (-1)^k \frac{N!}{k!} - \sum_{k=N+1}^\infty (-1)^k \frac{N!}{k!} = \frac{N!}{e} - R_N(-1).
\end{equation}
Here we used the notation $R_N(-1)$ for the remainder term of a Taylor expansion of $N! \, e^x$ around 0 at $x=-1$. Since the Taylor coefficients are simply $(-1)^k N!$, by the Remainder Estimation Theorem we have $|R_N(-1)| \leq (N+1)^{-1} \leq 1/3$ for all $N \geq 2$. Accordingly, \eqref{eq:FN-Euler} is always the integer closest to $N!/e$, or by using the floor function,
\begin{equation}\label{FN-short}
|F_N| = \left\lfloor \frac{N!}{e} + \frac{1}{2} \right\rfloor, \quad N \geq 2.
\end{equation}

\section{Fixed-point-free permutations having at least one $2$-cycle}

In this section we will only address permutations with at least one 2-cycle (one pair or more), which yields already part of the answer. We present three ways to calculate the number of fixed-point-free permutations with at least one 2-cycle, $P_N \subseteq S_N$, providing explicit formulas. 

\subsection{Cycle-types of permutations}

First, let us define the set of possible types of a permutation,
\begin{equation}
T_N = \left\{ (a_1,\ldots,a_N) \,\middle|\, a_i\geq 0, \sum\nolimits_{i=1}^N ia_i=N \right\}.
\end{equation}
We say that a permutation $\sigma \in S_N$ is of type $a=(a_1,\ldots,a_N)\in T_N$ if $\sigma$ consists of $a_i$ cycles of length $i$ (cycles of length $i=1$ just describe fixed-points, $i=2$ are the desired pairs).

The number of permutations of a given type can be calculated easily: For any permutation of type $a \in T_N$ we have $N$ open positions to fill with the elements $\{1,\dotsc,N\}$ where the slots are structured according to the given type. If for example $N=7$ and $a=(0,2,1,0,0,0,0)$ then the slots are like $(**)(**)(*{*}*)$. Now there are $N!$ ways how to interchange the slots, but many such variants will describe the same permutation and we have to get rid of those ambiguities. First, one can commute the $a_1$ cycles of length 1, the $a_2$ cycles of length 2 and so on, giving $a_1!\cdot\ldots\cdot a_N!$ such possibilities. Then inside each of the $a_i$ cycles of length $i$ one can write any of the $i$ elements in the beginning (canonical cycle notation with the largest element in front is not employed here), giving again $1^{a_1}\cdot\ldots\cdot N^{a_N}$ different possibilities of equivalent notation.
Hence, the number of permutations of type $a$ is $N! / (1^{a_1}\cdot\ldots\cdot N^{a_N} \cdot a_1!\cdot\ldots\cdot a_N!)$.

To answer our initial question, we get a first formula by summing over all permutations having zero $1$-cycles and at least one $2$-cycle: The size of the set of these permutations, which we denote by $P_N$, is thus given by
\begin{equation}\label{PN-sum}
|P_N|=\sum_{\substack{a\in T_N\\a_1=0\\a_2\geq 1}} \frac{N!}{1^{a_1}\cdot\ldots\cdot N^{a_N} \cdot a_1!\cdot\ldots\cdot a_N!}.
\end{equation}

Alas, for practical purposes this formula is complicated since the problem is just transferred into finding all the types considered in the sum.
To find out how many there are, we note that since every type considered in the sum in \eqref{PN-sum} fulfills $N=\sum_{i=2}^N i a_i$ (taking $a_1=0$ into account), this also gives an (increasing) partition of $N$ into integers $\geq 2$. According to the previous example we would thus have the partition $7 = 2+2+3$. Taking also $a_2 \geq 1$ into account such a partition of $N$ must include at least one $2$, so we start by considering the possible partitions of $N-2$ instead. All constituents of the partition of $N-2$ have to be greater or equal to $2$. For $2\leq v_1 \leq \dotsb \leq v_{k-1} \leq N-2$ let $(v_1,\dotsc,v_{k-1})$ be such a partition of $N-2$ into exactly $k-1$ numbers (one has already been accounted for). Mapping $(v_1,\dotsc,v_{k-1})\mapsto (v_1-2,\dotsc,v_{k-1}-2)$ delivers a partition of $N-2-2(k-1)=N-2k$ into \textit{at most} $k-1$ constituents. On the other hand, having a partition of $N-2k$ into at most $k-1$ constituents, $(u_1,\dotsc,u_{k-1})$ with $0\leq u_1\leq \dotsb\leq u_{k-1} \leq N-2k$, by mapping $(u_1,\dotsc,u_{k-1})\mapsto (u_1+2,\dotsc,u_{k-1}+2)$ we obtain a partition of $N-2$ into exactly $k-1$ constituents. This one-to-one relation between partitions of $N-2$ into exactly $k-1$ constituents that are greater or equal to $2$ and the partitions of $N-2k$ into at most $k-1$ constituents enables us to count them. In the example $2+2+3$ is mapped to a partition of $7-2\cdot 3=1$ which is arguably simple and tells us that there is only a single allowed partition of $7$ into three constituents greater or equal than $2$.
The possible number $k$ of constituents in any partition that meets our requirements can never exceed $\lfloor\frac{N}{2}\rfloor$, since each single constituent is $\geq 2$. This whole maneuver is meaningful because the number of partitions of $N$ into \emph{at most} $k$ constituents is given by the partition function\cite{wiki-partition-function} $p_k(N)$. So finally the number of summands in \eqref{PN-sum} is given by the formula
\begin{equation}
\sum_{k=2}^{\lfloor\frac{N}{2}\rfloor}p_{k-1}(N-2k)+1.
\end{equation}

Leaving the investigation of cycle types of permutations, next we look for another way to find an expression for our quantity $|P_N|$ which is hopefully easier to compute.

\subsection{A recursion for $|P_N|$}

We can find a way to enumerate the elements of $P_N$ by simply looking at a single number, say $N$. There are various, mutually exclusive possibilities for $N$: Either $N$ is part of a $2$-cycle, or $N$ is part of a $k$-cycle for exactly one $k$ with $3\leq k\leq N-2$. Note here, that no $N-1$ cycle is allowed, otherwise there would be a fixed-point, which we don't want. Concentrating on the first possibility, if $N$ in our permutation in $P_N$ is part of a $2$-cycle, then there are $N-1$ possibilities for the other element of this $2$-cycle. The remaining $N-2$ elements may form an arbitrary permutation without a fixed-point, which can be done in $|F_{N-2}|$ ways. So, there must be $|F_{N-2}|\cdot(N-1)$ permutations in $P_N$ such that $N$ is in a $2$-cycle. In the other circumstances, where $N$ belongs to some $k$-cycle, $3\leq k \leq N-2$, there are $\binom{N-1}{k-1}$ possibilities to choose the remaining $k-1$ elements of the cycle, forming $(k-1)!$ different cycles. The rest of the same permutation apart from this cycle must again be fixed-point-free and have at least one $2$-cycle, giving us $|P_{N-k}|$ possibilities. Summing up all those possibilities (which exhaust the set $P_N$) we get that for $N\geq 2$,
\begin{equation}\label{rec1}
\begin{aligned}
|P_N|&=\sum_{k=3}^{N-2} \binom{N-1}{k-1}(k-1)! \cdot |P_{N-k}| +|F_{N-2}|\cdot(N-1)\\
&=\sum_{k=2}^{N-3} |P_k|\cdot\frac{(N-1)!}{k!} +|F_{N-2}|\cdot(N-1),
\end{aligned}
\end{equation}
which is the first recursion for $|P_N|$, since we know that $|P_0|=|P_1|=0$. Using this recursion to calculate $|P_{N+1}|$ we immediately get that $|P_2|=1$ (which is what we expected) and, for $N\geq 4$,
\begin{align}
|P_{N}|&=(N-1)\sum_{k=2}^{N-4} |P_k|\frac{(N-2)!}{k!} +|P_{N-3}|(N-1)(N-2) +|F_{N-2}|\cdot(N-1)\\
&=(N-1)\left(|P_{N-1}|-|F_{N-3}|\cdot(N-2)\right)+|P_{N-3}|\cdot (N-1)(N-2)+|F_{N-2}|\cdot(N-1).\nonumber
\end{align}
We calculate the expression $|F_{N-2}|\cdot(N-1)-|F_{N-3}|\cdot(N-2)(N-1)$ to get
\begin{align}
    (N-1)!\sum_{k=0}^{N-2}\frac{(-1)^k}{k!}-(N-1)!\sum_{k=0}^{N-3}\frac{(-1)^k}{k!}=(N-1!)\frac{(-1)^{N-2}}{(N-2)!}=(-1)^N(N-1),
\end{align}
and hereby obtain the recursion for $|P_N|$,
\begin{align}\label{rec2}
    |P_{N}|=(N-1) \left( |P_{N-1}|+|P_{N-3}|\cdot (N-2)+(-1)^N \right),
\end{align}
which already seems quite practical to calculate some values of $|P_N|$:

\begin{center}
\begin{tabular}{r||c|c|c|c|c|c|c|c|c|c|c}
 $N$ & 0 & 1 & 2 & 3 & 4 & 5 & 6 & 7 & 8 & 9 & 10 \\ \hline
 $|P_N|$ & 0 & 0 & 1 & 0 & 3 & 20 & 105 & 714 & 5845 & 52632 & 525105
\end{tabular}
\end{center}

Although formula \eqref{rec2} looks simpler, in order to find a compact expression for $|P_N|$, we will use the seemingly more cumbersome form \eqref{rec1}. Note that because of $|P_0|=|P_1|=0$ we can write the first equality of \eqref{rec1} as
\begin{equation}\label{rec3}
|P_N|=\sum_{k=0}^{N} \frac{(N-1)!}{(N-k)!} \cdot |P_{N-k}| \cdot \chi_{k\geq 3} +|F_{N-2}|\cdot(N-1),
\end{equation}
where $\chi_{k\geq 3}=1$ if $k\geq 3$ and $0$ if $k\in \{0,1,2\}$. In the following, we will now use the technique of generating functions of sequences, where sequences $(f_n)_{n\geq 0}$ are represented by their according formal power series $\sum_{n=0}^\infty f_nx^n$ (see e.g.\ Chapter 3 in the book of M.~B\'ona \cite{bona}).
In order to apply this machinery, we multiply both sides of the recursion \eqref{rec3} with $x^N/(N-1)!$ and sum over $N=2,\ldots,\infty$ to get
\begin{equation}
    \sum_{N=2}^\infty |P_N| \cdot \frac{x^N}{(N-1)!}=\sum_{N=2}^\infty\sum_{k=0}^{N} \frac{x^N}{(N-k)!} \cdot |P_{N-k}| \cdot \chi_{k\geq 3} +\sum_{N=2}^\infty|F_{N-2}| \cdot (N-1) \frac{x^N}{(N-1)!}.
\end{equation}
Recognizing that the first sum on the right hand side can be started from $N=0$ because of $|P_0|=|P_1|=0$ and then can be rewritten as a Cauchy product, we find
\begin{equation}
\begin{aligned}
    x\sum_{N=1}^\infty |P_N|\cdot \frac{x^{N-1}}{(N-1)!}&=\left(\sum_{N=0}^\infty|P_N|\frac{x^N}{N!}\right) \cdot \left(\sum_{N=0}^{\infty}\chi_{N\geq 3}\cdot x^N\right)+\sum_{N=2}^\infty|F_{N-2}|\frac{x^N}{(N-2)!}\\
    &=\left(\sum_{N=0}^\infty|P_N|\frac{x^N}{N!}\right)\cdot\frac{x^3}{1-x}+x^2\sum_{N=0}^\infty|F_{N}|\frac{x^N}{N!}\\
    &=\left(\sum_{N=0}^\infty|P_N|\frac{x^N}{N!}\right)\cdot\frac{x^3}{1-x}+x^2\left(\sum_{N=0}^\infty\sum_{k=0}^N\frac{(-1)^k}{k!}x^N\right)\\
    &=\left(\sum_{N=0}^\infty|P_N|\frac{x^N}{N!}\right)\cdot\frac{x^3}{1-x}+x^2\left(\sum_{N=0}^\infty\frac{(-x)^N}{N!}\right)\cdot\left(\sum_{N=0}^\infty x^N\right)\\
    &=\left(\sum_{N=0}^\infty|P_N|\frac{x^N}{N!}\right)\cdot\frac{x^3}{1-x}+\frac{x^2e^{-x}}{1-x}.
\end{aligned}
\end{equation}
We cancel $x$ and denote $p(x)=\sum_{N=0}^\infty|P_N|x^N/N!$, which leads us to the ordinary differential equation
\begin{equation}
\frac{\mathrm{d}}{\mathrm{d} x}p(x)=p(x)\cdot \frac{x^2}{1-x}+\frac{x e^{-x}}{1-x},\quad p(0)=0,
\end{equation}
which can be solved by the formula of variation of constants, using the integral $\int_0^x y^2/(1-y) \,\mathrm{d}y=-x^2/2-x-\log(1-x)$ that holds for $x< 1$, to end up with
\begin{equation}\label{eq:gf-no-fixed-points}
\begin{aligned}
    p(x)&=\exp\left(-\frac{x^2}{2}-x-\log(1-x)\right) \int_0^x \frac{ye^{-y}}{1-y}\cdot \exp\left(\frac{y^2}{2}+y+\log(1-y)\right)\,\mathrm{d} y\\
    &=\frac{e^{-x^2/2-x}}{1-x}\int_0^x y\cdot e^{y^2/2} \,\mathrm{d} y=\frac{e^{-x}}{1-x}-\frac{e^{-x^2/2-x}}{1-x}.
\end{aligned}
\end{equation}
The series expansion of the resulting function above can be written (multiplying out the Cauchy products) as
\begin{equation}
    \sum_{N=0}^\infty|P_N|\frac{x^N}{N!}=p(x)=\sum_{N=0}^\infty \sum_{k=0}^N\frac{(-1)^k}{k!}x^N-\sum_{N=0}^\infty \sum_{k=0}^N\sum_{j=0}^k \frac{(-1)^{k-j}(-1)^{j/2}}{(k-j)!\,\frac{j}{2}!\,2^{j/2}}\chi_{2\mathbb{N}}(j) \cdot x^N.
\end{equation}
Here $\chi_{2\mathbb{N}}$ denotes the characteristic function of the even numbers. From this we obtain the explicit formula for $|P_N|$ term by term
\begin{align}\label{rec4}
    |P_N|=N!\sum_{k=0}^N\frac{(-1)^k}{k!}-N!\sum_{k=0}^N \frac{(-1)^{k-j}(-1)^{j/2}}{(k-j)!\,\frac{j}{2}!\,2^{j/2}}\chi_{2\mathbb{N}}(j).
\end{align}
It is no mere coincidence that the first sum here is $|F_N|$ as we will also see in the next section, where we achieve the same result in \eqref{PN-res} with a different method.

\subsection{A further method to obtain $|P_N|$ and asymptotic results}

The formula for $|P_N|$ found in the end of the last paragraph is explicit and ready to use for further mathematical treatment, which we will do in this section. Before however, we will use a theorem from enumerative combinatorics which leads to a quick way to compute the expression for $|P_N|$. The key idea in this section is to look at the set of permutations that have no fixed-points and neither possess $2$-cycles. The cardinality of this set can be obtained via Theorem 4.34 from the combinatorics book of M.~B\'ona\cite{bona} (and equals the second summand of \eqref{rec4}).

\begin{theorem*}
Let $M$ be any set of positive integers and let $g_M(n)$ be the number of permutations of length $n$ whose cycle lengths are all elements of $M$. Then
\begin{equation}\label{bona-thm-eq}
G_M(x)=\sum_{N=0}^\infty g_M(N)\frac{x^N}{N!}=\exp\left(\sum_{k\in M}\frac{x^k}{k!}\right).
\end{equation}
\end{theorem*}

In the case of permutations having cycle lengths strictly larger than $2$, the set $M$ must be $\{3,4,\dotsc\}$. Hence, 
\begin{equation}
G_{\geq 3}(x)=\exp\left(\sum_{k=3}^\infty \frac{x^k}{k!}\right),
\end{equation}
by the above theorem. With the Taylor-series expansion $-\log(1-x)=\sum_{k=1}^\infty\frac{x^k}{k!}$ that holds for $-1 \leq x < 1$, we get
\begin{equation}
G_{\geq 3}(x)=\exp\left(-\log(1-x)-x-\frac{x^2}{2}\right)=\frac{\exp\left(-x-\frac{x^2}{2}\right)}{1-x}.
\end{equation}
We now use the series expansion of the exponential and the geometric series to arrive at
\begin{equation}
G_{\geq 3}(x)=\left(\sum_{k=0}^\infty x^k\right)\left(\sum_{k=0}^\infty(-1)^k\frac{x^k}{k!}\right)\left(\sum_{k=0}^\infty\frac{(-1)^k x^{2k}}{2^k k!}\right).
\end{equation}
Multiplying out the appearing Cauchy products and simplifying the arising sums we get that the coefficient $g_{\geq 3}(N)$ from \eqref{bona-thm-eq} equals
\begin{equation}\label{formula}
g_{\geq 3}(N) = N!\sum_{k=0}^N \sum_{j=0}^k \frac{(-1)^{k-j}(-1)^{j/2}}{(k-j)! \, \frac{j}{2}!\,2^{j/2}}\chi_{2\mathbb{N}}(j).
\end{equation}
The inner sum can be further contracted to
\begin{equation}
g_{\geq 3}(N)=N!\sum_{k=0}^N (-1)^k \sum_{j=0}^{\lfloor\frac{k}{2}\rfloor}\frac{(-1)^j}{(k-2j)!\,j!\,2^j}.
\end{equation}
To find the number $|P_N|$ we now just have to subtract the obtained quantity from the number of all permutations without fixed-points \eqref{FN-num}, so
\begin{equation}\label{PN-res}
|P_N|=|F_N|-g_{\geq 3}(N).
\end{equation}

\begin{figure}
  \includegraphics[width=.8\linewidth]{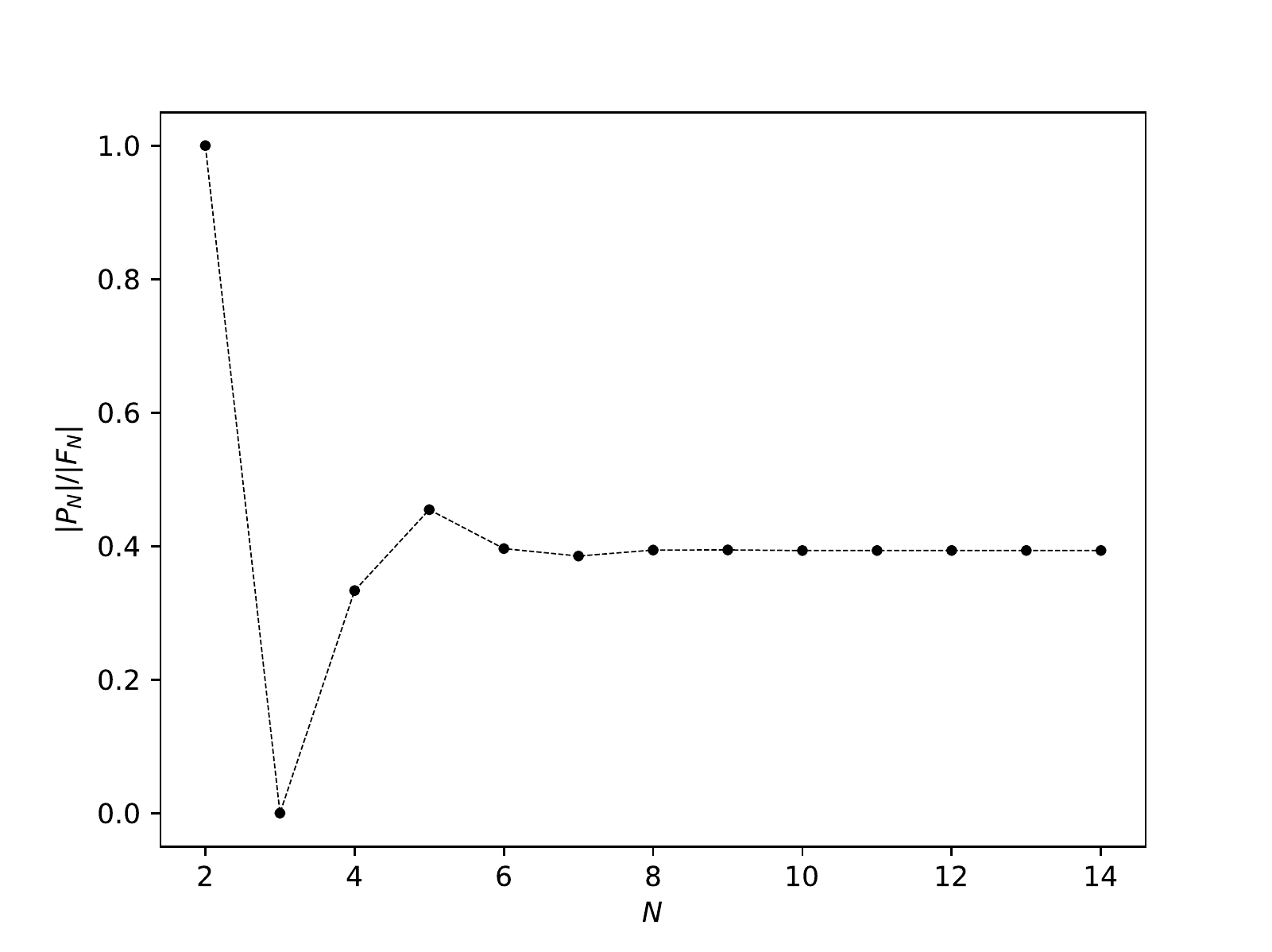}
  \captionof{figure}{probability for at least one pair}
  \label{fig1}
\end{figure}

Having obtained a decent formula, we are ready to calculate the asymptotic probability to have at least one pair: We know that there are overall $|F_N|$ possible Secret Santa configurations. The relative part of configurations where at least one pair exchanges gifts is then $|P_N|/|F_N|$ which equals
\begin{equation}\label{PN-FN}
\frac{|P_N|}{|F_N|}=1 - \frac{g_{\geq 3}(N)}{|F_N|}.
\end{equation}
In the limit $N \to \infty$ the two sums of \eqref{formula} can be split into a Cauchy product again and one has
\begin{equation}
\lim_{N \to \infty} \frac{g_{\geq 3}(N)}{N!} = \left( \sum_{k=0}^\infty \frac{(-1)^k}{k!} \right) \left( \sum_{k=0}^\infty \frac{(-1)^{k/2}}{\frac{k}{2}!\,2^{k/2}}\chi_{2\mathbb{N}}(k) \right).
\end{equation}
That means in the asymptotics of \eqref{PN-FN} the $|F_N|$ given by \eqref{FN-num} cancels and we are left with
\begin{equation}\label{asymp-pair}
\lim_{N \to \infty} \frac{|P_N|}{|F_N|} = 1- \sum_{k=0}^\infty \frac{(-1)^{k/2}}{\frac{k}{2}!\,2^{k/2}}\chi_{2\mathbb{N}}(k) = 1-e^{-1/2}\approx 0.393.
\end{equation}
If the whole world plays Secret Santa, we can expect at least one gift-exchanging with a probability of about 40\%. As one can see in Fig.~\ref{fig1} convergence is really fast.

\section{Fixed-point-free permutations with exactly $k$ $2$-cycles}

Going back to the initial question, we are not only interested in the number of configurations in Secret Santa that admit at least one pair, but in the number where exactly $k$, $0\leq k\leq \lfloor \frac{n}{2} \rfloor$, pairs exchange gifts. The set of fixed-point-free permutations of $N$ elements that have exactly $k$ $2$-cycles will be called $P_N^k$. This will lead us then to a final answer. Now how big are the probabilities that such things happen?
Let us tackle this question with a similar strategy as when we investigated permutations of a given cycle type. We can think about a permutation being split up into $k$ $2$-cycles and the remaining ($2$-cycle-free) permutation. In choosing the first element of the first $2$-cycle there are $N$ possibilities, for the second element of the first $2$-cycle there are $N-1$ possibilities, for the first element of the second $2$-cycle then $N-2$ possibilities and so on, until we have $(N-(2k-1))$ possibilities left for the second element of the $k$-th $2$-cycle. In each one of the $k$ $2$-cycles any of its two elements can be written first, so there are $2^k$ ways to describe the same $2$-cycles. Moreover, due to the commutativity of cycles there are $k!$ ways that the $2$-cycles just constructed can be arranged equivalently. Summarizing, we have
\begin{equation}
\frac{N(N-1)\dotsm(N-(2k-1))}{2^k k!}=\frac{N!}{(N-2k)!2^k k!}
\end{equation}
ways to form $k$ $2$-cycles with $2k$ of the $N$ numbers. For the remaining permutation of $N-2k$ elements without $1$- or $2$-cycles, there are $g_{\geq 3}(N-2k)$ possibilities, with $g_{\geq 3}$ from \eqref{formula}. Hence, the number of fixed-point-free permutations with exactly $k$ $2$-cycles is
\begin{equation}\label{PNk-res}
|P_N^k| = \frac{N!}{(N-2k)!2^k k!}\cdot g_{\geq 3}(N-2k)=\frac{N!}{2^k k!}\sum_{m=0}^{N-2k} \sum_{j=0}^m \frac{(-1)^{m-j}(-1)^{j/2}}{(m-j)! \, \frac{j}{2}!\,2^{j/2}}\chi_{2\mathbb{N}}(j).
\end{equation}
In the above formula we can correctly check that $|P_N^0| = g_{\geq 3}(N)$, the number of fixed-point-free permutations with no (exactly zero) $2$-cycles.
Hence, the probability to have exactly $k$ pairs that exchange presents
is $|P_N^k|/|F_N|$. With the same type of manipulations that lead to \eqref{asymp-pair} we arrive at the asymptotic probability
\begin{equation}
\lim_{N \to \infty}\frac{|P_N^k|}{|F_N|} = \frac{\left(\frac{1}{2}\right)^k e^{-1/2}}{k!}.
\end{equation}
This shows that, asymptotically, the probability that exactly $k$ pairs draw each other follows a Poisson distribution with parameter $\frac{1}{2}$. In the very case that $3$ pairs of people exchange gifts when playing Secret Santa with a high number of persons, we get a small probability of approximately $0.013$.

We can (via e.g.\ the Skorohod representation\cite{skorohod-representation}) assume that there are (discrete) random variables $(X_N)_{N\geq 2}$ (the question is not very interesting for $N = 0$ or $1$) converging almost surely to a random variable $X$, where $\mathbb{P}(X_N = k)=|P_N^k|/|F_N|$ and $\mathbb{P}(X = k)= e^{-1/2}/(2^k k!)$. Thus the random variable $X_N$ describes the number of pairs found playing Secret Santa with $N$ players and $X$ describes the same for the limiting case (which would correspond to playing it with infinitely many players). Let us estimate for $s\in\mathbb{R}$ the expectation
\begin{align}
    \mathbb{E}\left( e^{s X_N} \right)&=\sum_{k=0}^\infty e^{sk}\frac{|P_N^k|}{|F_N|}=\sum_{k=0}^\infty \frac{e^{sk}}{2^k k! \sum_{m=0}^N \frac{(-1)^m}{m!}}\sum_{m=0}^{N-2k} \sum_{j=0}^m \frac{(-1)^{m-j}(-1)^{j/2}}{(m-j)! \, \frac{j}{2}!\,2^{j/2}}\chi_{2\mathbb{N}}(j)\nonumber \\
    & \leq \sum_{k=0}^\infty \frac{e^{sk}}{2^k k! \frac{1}{3}} \sum_{m=0}^{\infty}\sum_{j=0}^m \frac{\chi_{2\mathbb{N}}(j)}{(m-j)! \, \frac{j}{2}!\,2^{j/2}}=3\sum_{k=0}^\infty \frac{e^{sk}}{2^k k!} \sum_{m=0}^{\infty}\frac{1}{m!}\sum_{m=0}^\infty \frac{1}{2^m m!}\nonumber \\
    &=  3\sum_{k=0}^\infty \frac{e^{sk} \cdot e \cdot e^\frac{1}{2}}{2^k k!} =3e^{e^{s/2}+\frac{3}{2}},
\end{align}
where for the inequality we used that the double sum can be estimated by the infinite series adding up the absolute value of the summands. We further used that the Leibniz partial sums $\sum_{m=0}^N \frac{(-1)^m}{m!}$ for $N\geq 2$ are always larger than the one for the smallest odd $N\geq 2$, which is $1-1+\frac{1}{2}-\frac{1}{6}=\frac{1}{3}$.
Therefore, $\sup_{N\geq 2} \mathbb{E}\left( e^{s X_N}\right) <\infty$ for all $s\in \mathbb{R}$ and by dominated convergence we get that all (even the exponential) moments of $X_N$ converge to the moments of $X$ that are those of a Poisson distribution with parameter $\tfrac{1}{2}$. 
The expected number of $2$-cycles in the long run is then just the parameter of the Poisson distribution, which is $\frac{1}{2}$.

The final answer is still due: How large is the probability for three \emph{or more} pairs? In this case we just have to eliminate the configurations with exactly zero, one, or two pairs. This is now an easy task applying the above result, and we get
\begin{equation}\label{final-prob}
\lim_{N \to \infty} \frac{|F_N| - |P_N^0| - |P_N^1| - |P_N^2|}{|F_N|} = 1 - \left(1 + \frac{1}{2} + \frac{1}{8} \right) e^{-1/2} = 1 - \frac{13 \, e^{-1/2}}{8} \approx 0.014
\end{equation}

In the case of the choir $N=20$ persons and the probability from \eqref{final-prob} is already equal to the asymptotic one up to 7 digits after the comma. This fast convergence rate is visible in Fig.~\ref{fig2}, where also the peak of probability for $N=6$ is revealed.

\begin{center}
  \includegraphics[width=.75\linewidth]{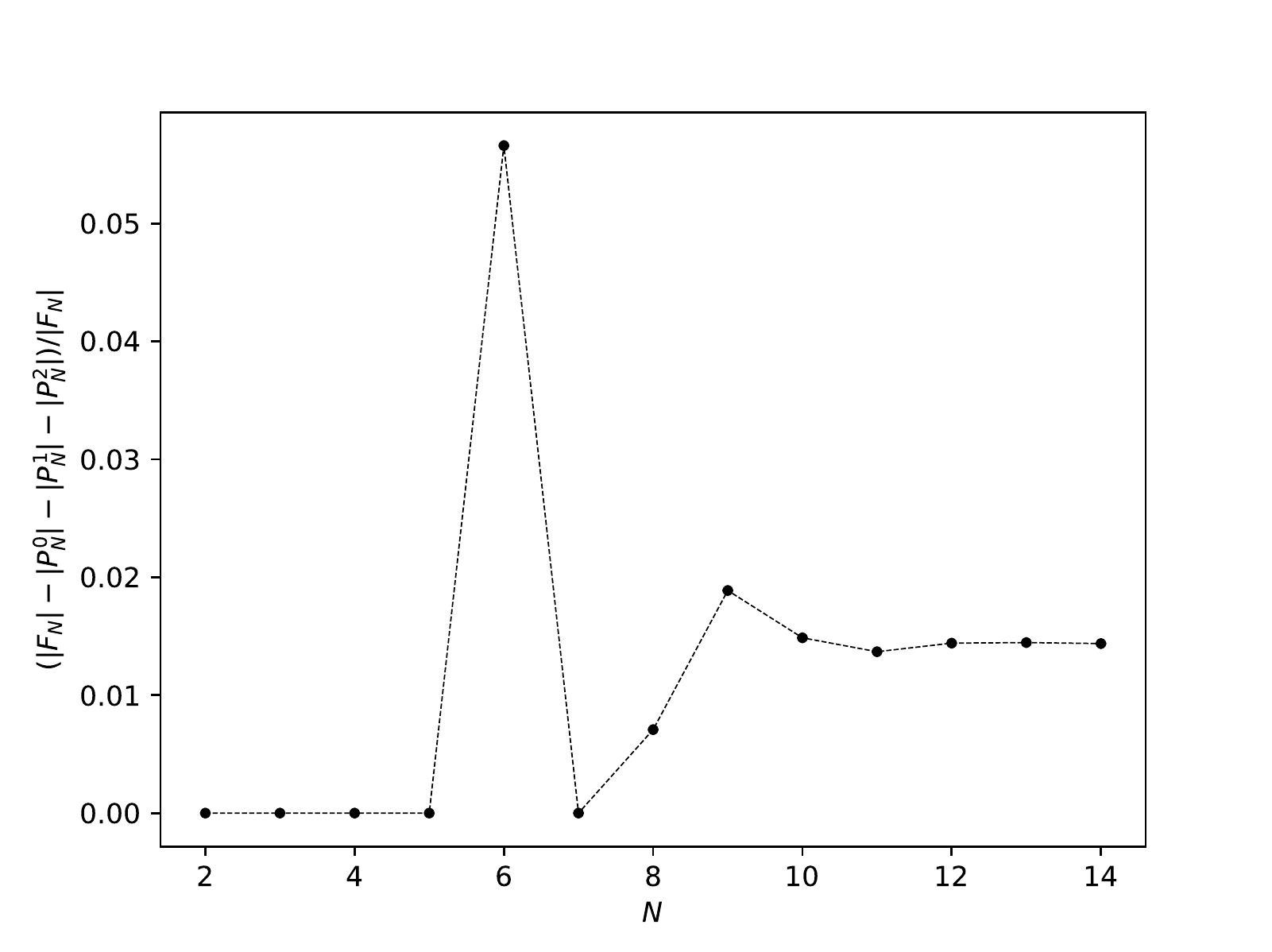}
  \captionof{figure}{probability for at least three pairs}
  \label{fig2}
\end{center}

\end{document}